\documentclass{ws-procs9x6}
\usepackage{mathrsfs}

\usepackage[all]{xy}
\usepackage{amssymb}
\usepackage{amsmath}

\def\Ai{A_{\infty}}

\def\d{\succ}
\def\g{\prec}

\def\DD{\Delta}

\def\t{\otimes}

\def\RR{{\mathbb{R}}}
\def\CC{{\mathbb{C}}}

\def\HHH{{\mathbb{H}}}
\def\OOO{{\mathbb{O}}}

\def\CCC{{\mathcal{C}}}

\def\PP{{\mathcal{P}}}

\def\Id{\mathrm{Id }}

\def\Vect{\mathop{\rm Vect }}

\def\id{\mathrm{ id }}

\def\KK{\mathbb{K}}

\def\arbreA{\vcenter{\xymatrix@R=3pt@C=3pt{
&& \\
&*{}\ar@{-}[ur] \ar@{-}[ul] \ar@{-}[d]     &\\
&&
}}}

\begin{document}

\title{Some problems in operad theory}

\author{Jean-Louis Loday}
\address{Institut de Recherche Math\'ematique Avanc\'ee\\
    CNRS et Universit\'e de Strasbourg\\
    7 rue R. Descartes\\
    67084 Strasbourg Cedex, France\\
e-mail: loday@math.unistra.fr
}

\begin{abstract} This is a list of some problems and conjectures related to various types of algebras, that is to algebraic operads. Some comments and hints are included.
\end{abstract}


\bodymatter

\section*{Introduction} Since 1991 I got involved in the operad theory, namely after an enlightening lecture by Misha Kapranov in Strasbourg. During these two decades I came across many questions and problems. The following is an excerpt of this long list which might be helpful to have in mind while working in this theme. Of course this is a very personal choice.

Notation and terminology are those of [\cite{Lod08}] and [\cite{LV11}]. Various types of algebras, i.e.\ algebraic operads, can be found in [\cite{Zinb}].

\section{On the notion of group up to homotopy} The notion of associative algebra up to homotopy is well-known: it is called \emph{$\Ai$-algebra} and was devised by Jim Stasheff in [\cite{Sta63}]. It has the following important property: starting with a differential graded associative algebra $(A,d)$, if $(V,d)$ is a deformation retract of $(A,d)$, then $(V,d)$ is not a dg associative algebra in general, but it is an $\Ai$-algebra. This is Kadeishvili's theorem [\cite{Kad82}], see [\cite{LV11}] for a generalization and variations of it. It is called the Homotopy Transfer Theorem. Let us now start with a group $G$. What is the notion of a ``group up to homotopy''~? To make this question more precise we move, as in quantum group theory, to the group algebra $\KK[G]$ of the group $G$ over a field $\KK$. It is well-known that this is not only a unital associative algebra, but it is a cocommutative Hopf algebra. So, it has a cocommutative coproduct $\DD$ (induced by the diagonal on $G$), and the existence of an inverse in $G$ translates to the existence of an antipode on the group algebra. So we can now reformulate the question as follows:

\medskip

``What is the notion of cocommutative Hopf algebra up to homotopy ?''

\medskip

One of the criterions for the answer to be useful would be the existence of a Homotopy Transfer Theorem for cocommutative Hopf algebras. One has to be careful enough to take into account that the existence of a unit and a counit is part of the structure of a Hopf algebra. In the associative case the operad $\Ai$ does not take the unit into account. See [\cite{Lod08}] for the Hopf relation of a nonunital bialgebra.

The fact that the tensor product of two associative algebras is still an associative algebra plays a prominent role in the definition of a bialgebra (a fortiori a Hopf algebra). So it is clear that a first step in analyzing this problem is to check whether one can put an $\Ai$-structure on the tensor product of two $\Ai$-algebras and to unravel the properties of such a construction. A first answer has been given by Saneblidze and Umble in [\cite{SU}]. But this tensor product is not associative. This problem has been addressed in [\cite{Lod11}].

\section{Subgroup of free group} It is well-known that a subgroup of a free group is free. The proof is topological in the sense that it consists in letting the free group act on a tree. Could one find a proof by looking at the properties of the associated group algebra (which is a Hopf algebra) ? I am thinking about something similar to the theorems which claim that some algebra is free under certain condition (PBW type theorems, see [\cite{Lod08}]).

\section{The octonions as an algebra over a Koszul operad} The octonions form a normed division algebra $\OOO$ of dimension $8$, see for instance [\cite{Baez02}]. The product is known not to be associative contrarily to the other normed division algebras $\RR, \CC$ and $\HHH$. However it does satisfy some algebraic relation: it is an alternative algebra. Let us recall that an alternative algebra is a vector space equipped with a binary operation $x\cdot y$, such that the associator $(x,y,z):=(x\cdot y)\cdot z - x\cdot (y\cdot z) $ is antisymmetric:
$$(x,y,z)=-(y,x,z)=-(x,z,y).$$

It turns out that the operad of alternative algebras is not too good, because it is not a Koszul operad, cf.\ [\cite{DZ}]. Whence the question: Find a (small) binary operad such that the octonions form an algebra over this operad and such that this operad is Koszul. Of course the ambiguity of the question is in the adjective ``small''. Because such an operad exists: it suffices to take the magmatic operad on one binary operation. But there may exists a smaller operad (i.e.\ a quotient of $Mag$), which is best. This operad need not be quadratic, that is, we may look for relations involving 4 variables, like in Jordan algebras.

\section{Commutative algebras up to homotopy in positive characteristic} In positive characteristic $p$ it is best to work with divided power algebras rather than commutative algebras. The notion of commutative algebra up to homotopy is well-known: it is the $C_{\infty}$-algebras (also denoted $Com_{\infty}$, see [\cite{LV11}]). What is, explicitly, the notion of divided power algebra up to homotopy in characteristic $p$ ? Theoretically the problem can be solved as follows. One can perform the theory of Koszul duality for operads with divided powers, cf.\ [\cite{Fre00}]. Any such  operad with divided powers $\Gamma \PP$, which is Koszul, gives rise to a dg  operad with divided powers $\Gamma \PP_{\infty}$. A divided power algebra up to homotopy is an algebra over $\Gamma \PP_{\infty}$. The point is to make all the steps of the theory explicit in the case $\PP= Com$.

\section{Manin black product for operads} What is the operad $Com \bullet Ass$ ? One is asking for a small presentation by generators and relations.

\section{$L$-dendriform algebras and operadic black product} By definition an $L$-dendriform algebra is a vector space equipped with 
two operations $x\g y$ and $x\d y$ satisfying 
\begin{eqnarray*}(x \g y) \g  z + y \d  (x \g  z) &=& x \g  (y * z) + (y \d  x) \g  z ,\\
 (x * y) \d  z + y \d  (x \d  z) &=&x \d  (y \d  z) + (y \g x) \d  z ,
\end{eqnarray*}
where $x * y = x \g  y + x \d  y$ (cf.\ [\cite{BLN}]).  This is one of the numerous ways of splitting the associativity of the operation $*$. 

\bigskip

\noindent {\bf Conjecture}:  the operad encoding $L$-dendriform algebras is a Manin black product:
$$preLie \bullet preLie = L\textrm{-}Dend.$$
In favor of this conjecture we have the following facts (cf.\ [\cite{Val08}]):
$$preLie \bullet Com = Zinb,\  preLie \bullet Ass = Dend,\ preLie \bullet Lie = preLie,$$
and also $preLie \bullet Dend = Quad$. 

Similarly it seems that $preLie \bullet Zinb = ComQuad$ and that the operads $Octo, L\textrm{-}Quadri, L\textrm{-}Octo$ introduced in [\cite{LNB}] are also black products:
\begin{gather*}
preLie \bullet Quadri = Octo,\\
preLie \bullet L\textrm{-}Dend=L\textrm{-}Quadri,\\
 preLie \bullet L\textrm{-}Quadri=L\textrm{-}Octo.
 \end{gather*}

\noindent{\bf Note.} Some of these questions have been recently settled in [\cite{BBGN}].

\section{Resolutions of associative algebras} Koszul duality theory for associative algebras gives a tool to construct free resolutions for some associative algebras, and even the minimal resolutions in certain cases. When the algebra is a group algebra, then there are tools to construct (at least the beginning of) a resolution by taking the free module on the set of generators, then of relations, then of relations between the relations, and so forth (syzygies), see for instance [\cite{Lod00}]. It would be very interesting to compare these various methods.

\section{Hidden structure for EZ-AW maps} Given a deformation retract 
\begin{eqnarray*}
&\xymatrix{     *{ \quad \ \  \quad (A, d_{A})\ \ } \ar@(dl,ul)[]^{h}\ \ar@<1ex>[r]^{p} & *{\quad
(V,d_{V})\quad \ \  \ \quad } \ar@<1ex>[l]^{i}},&\\
pi = \id_{V},&\Id_A-ip =d_{A}  h + h  d_{A},&
\end{eqnarray*}the HTT says that an algebraic structure on $(A,d_{A})$ can be transferred to some other algebraic structure (the hidden one) on $(V,d_{V})$. This principle is not special to chain complexes and can be applied to other situations as shown for crossed modules in [\cite{LV11}]. Apply this principle to the Eilenberg-Zilber and Alexander-Whitney quasi-isomorphisms. Let us  recall that, for $X$ and $Y$ being simplicial modules,  these isomorphisms relate the chain complex $C_{\bullet}(X\times Y)$ to the tensor product $C_{\bullet}(X)\t C_{\bullet}(Y)$. 

\section{Interpolating between $Dend$ and $Com$} The so-called $E_{n}$-operads are operads which interpolate between the homotopy class of the operad $Ass$, which contains $\Ai$, and the homotopy class of the operad $Com$, which contains $C_{\infty}$. So, it solves the question: what is an associative algebra which is more or less commutative. The answer is: an $E_{n}$-algebra ; the larger $n$ is,  the more commutative it is. 

Question: what is a dendriform algebra which is more or less commutative ? In other words we are looking for an interpolation $D_{n}$ between the operads $Dend_{\infty}$ and $Zinb_{\infty}$, where $Dend$ is the operad of dendriform algebras (two generating operations $\g$ and $\d$ and three relations) and $Zinb$ is the operad of Zinbiel algebras (dendriform algebras such that $x\d y= y\g x$).

One of the motivation for finding $D_{2}$ is the following. It is known that the Grothendieck-Teichm\" uller group is related to the operad $E_{2}$. Knowing $D_{2}$ could lead to a dendriform version of the 
Grothendieck-Teichm\" uller group.

\section{Good triples of binary quadratic operads} Let $\PP$ be a binary quadratic operad which is Koszul (cf.\ for instance [\cite{MSS, LV11}]). It gives a notion of $\PP$-algebra and also a notion of $\PP$-coalgebra. We conjecture that there is a compatibility relation which defines a notion of $\PP^{c}\textrm{-}\PP$-bialgebra such that $(\PP, \PP, \Vect)$ is a good triple of operads in the sense of [\cite{Lod08}].

\noindent \texttt{Comments.} There are many examples known: 
$$\PP= Com, As, Dend, Mag, 2\textrm{-}as.$$
 When it holds, it gives a criterion for proving that a given $\PP$-algebra is free.

\section{On the coalgebra structure of Connes-Kreimer Hopf algebra}  The Connes-Kreimer Hopf algebra is an algebra of polynomials endowed with an ad hoc coproduct, cf.\ [\cite{CK}]. It is known that the indecomposable part is not only coLie, but in fact co-pre-Lie, cf. [\cite{CL}]. If we linearly dualize (as graded modules), the Hopf algebra is the Grossman-Larson Hopf algebra, which is cocommutative and its primitive part is pre-Lie. I conjecture that there is some type of algebras, that is some operad $\mathcal{X}$, and  some type of $Com^{c}\textrm{-} \mathcal{X}$-bialgebras, which fit into a good triple of operads
$$(Com, \mathcal{X}, preLie).$$
If so, then the Grossman-Larson algebra would be the free $\mathcal{X}$-algebra on one generator.

The solution of this problem in the noncommutative framework is given by the operad $Dend$, cf.\ [\cite{LR11}].

\section{Generalized bialgebras in positive characteristic} The Poincar\' e-Birkhoff-Witt theorem and the Cartier-Milnor-Moore theorem are structure theorems for cocommutative bialgebras in characteristic zero. They can be summarized by saying the triple of operads
$$(Com, As, Lie)$$
is a good triple. Several other good triples have been described in [\cite{Lod08}], some of them being valid in any characteristic, like the triple $(As, Dup, Mag)$ for instance. For the classical case, it is known that, in order for the CMM theorem to be true in characteristic $p$, one has to replace the notion of Lie algebra by the notion of \emph{restricted Lie algebras}. Operadically, restricted Lie algebras, divided power algebras and the like are obtained by replacing the ``coinvariants'' in the definition of an operad by the ``invariants'', cf. \ [\cite{Fre00}],
$$\Gamma\PP(V):= \sum_n (\PP(n)\t V^{\t n})^{S_n}.$$
So the PBW-CMM theorem in characteristic $p$ can be phrased by saying that 
$$(Com, As, \Gamma Lie)$$
is a good triple. Note that $As= \Gamma As$, so equivalently $(Com, \Gamma As, \Gamma Lie)$ is a good triple.

It would be very interesting to generalize the results on generalized bialgebras to positive characteristic along these lines, that is, to show that, when $(\CCC, \mathcal{A}, \PP)$ is a good triple, then so is $(\CCC, \Gamma \mathcal{A}, \Gamma \PP)$.

Similarly, $(\Gamma Com, As, Lie)$ is a good triple. One should be able to show that, when $(\CCC, \mathcal{A}, \PP)$ is a good triple, then so is $(\Gamma \CCC, \mathcal{A},\PP)$.

\section{Higher Dynkin diagrams and operads} Show that there exists some types of algebras (i.e.\ some operads) for which the finite dimensional simple algebras are classified by the diagrams described by Ocneanu in [\cite{Ocn02}]. The toy-model is the operad $Lie$ and the Dynkin diagrams. 

\section{Coquecigrues} It is well-known that a Lie group admits a tangent space at the unit element which is a Lie algebra. But there is also another relationship between groups and Lie algebras, more specifically between \emph{discrete groups} and Lie algebras (over $\mathbb{Z}$). It is given by the descending central series. For $G$ a discrete group, $G^{(1)}=[G,G]$ is its commutator subgroup, and, more generally, $G^{(n)}= [G, G^{(n-1)}]$ is the $n$th term of the descending central series. It is well-known that the graded abelian group $\bigoplus_{n}G^{(n)}/G^{(n+1)}$ is a Lie algebra whose bracket is induced by the commutator in $G$. The Jacobi identity is a consequence of a nice (and not so well-known) relation, valid in any group $G$, called the Philip Hall relation (see for instance [\cite{Lod00}] for some \emph{drawing} of it related to the Borromean rings). 

A natural question is the following. Let $\PP$ be a variation of the operad $Lie$ (we have in mind $preLie$ and $Leib$). Is there some structure playing the role of groups in this realm ? For Leibniz algebras the question arised naturally in my research on the periodicity properties of algebraic $K$-theory, cf. [\cite{Lod87, Lod03}]. I called this conjectural object a \emph{coquecigrue}. In fact I was more interested in the cohomology theory which should come with this new notion, to apply it further to groups. Recent progress using the notion of racks was achieved by Simon Covez in [\cite{Covez10}].

\section{Homotopy groups of spheres} Let $p$ be a prime number. Let $(\pi^{S}_{{\bullet}}(X))_{p}$ be the stable homotopy groups of the pointed connected topological space $X$, localized at $p$. Find a type of algebras such that $(\pi^{S}_{{\bullet}}(X))_{p}$ is an algebra of this type and such that $(\pi^{S}_{{\bullet}}(*))_{p}$ is the free algebra of this type over one generator (in degree $2p-3$). 

\bigskip

\noindent \texttt{Comments.} The Toda brackets are likely to play a role in this problem.

\bibliographystyle{alpha}

\end{document}